\documentclass{article}

\usepackage{amsmath,amssymb,amsthm}

\newtheorem{theorem}{Theorem}[section]
\newtheorem{definition}[theorem]{Definition}
\newtheorem{proposition}[theorem]{Proposition}
\newtheorem{corollary}[theorem]{Corollary}
\newtheorem{lemma}[theorem]{Lemma}
\newtheorem{fact}[theorem]{Remark}
\newtheorem{exemplu}[theorem]{Example}

\newtheorem{remark}[theorem]{Remark}
\newcommand{\bdfn}{\begin{definition}}
\newcommand{\edfn}{\end{definition}}
\newcommand{\bthm}{\begin{theorem}}
\newcommand{\ethm}{\end{theorem}}
\newcommand{\bprop}{\begin{proposition}}
\newcommand{\eprop}{\end{proposition}}
\newcommand{\bcor}{\begin{corollary}}
\newcommand{\ecor}{\end{corollary}}
\newcommand{\blem}{\begin{lemma}}
\newcommand{\elem}{\end{lemma}}
\newcommand{\bfact}{\begin{fact}}
\newcommand{\efact}{\end{fact}}
\newcommand{\bex}{\begin{exemplu}\begin{rm}}
\newcommand{\eex}{\end{rm}\end{exemplu}}

\def\R{{\mathbb R}}
\def\N{{\mathbb N}}

\def\Q{{\mathbb Q}}

\newcommand{\nin}{\!\in\!\!\!\!\!\!/\,}

\newcommand{\lambdaxy}{(1-\lambda)x\oplus\lambda y}

\newcommand{\eps}{\varepsilon}
\newcommand{\ol}{\overline}

\newcommand{\be}{\begin{enumerate}}
\newcommand{\ee}{\end{enumerate}}
\newcommand{\bt}{\begin{tabular}}
\newcommand{\et}{\end{tabular}}
\newcommand{\beq}{\begin{equation}}
\newcommand{\eeq}{\end{equation}}
\newcommand{\ba}{\begin{array}} 
\newcommand{\ea}{\end{array}}
\newcommand {\bea} {\begin{eqnarray}}
\newcommand {\eea} {\end {eqnarray}}
\newcommand {\bua} {\begin{eqnarray*}}
\newcommand {\eua} {\end {eqnarray*}}
\newcommand{\se}{\subseteq}
\newcommand{\ds}{\displaystyle}

\newcommand{\ra}{\rightarrow}
\newcommand{\si}{\wedge}

\begin{document}

\title{Asymptotically nonexpansive mappings in uniformly convex hyperbolic 
spaces\thanks{The research reported in this paper was carried out during 
the authors stay 
at the Max-Planck-Institute for Mathematics (Bonn) whose support is gratefully 
acknowledged.}}
\author{U. Kohlenbach$^{1}$, L. Leu\c stean$^{1,2}$\\[0.2cm]
\footnotesize ${}^1$ Department of Mathematics, Technische Universit\" at Darmstadt,\\
\footnotesize Schlossgartenstrasse 7, 64289 Darmstadt, Germany\\[0.1cm]
\footnotesize${}^2$ Institute of Mathematics "Simion Stoilow'' of the 
Romanian Academy, \\
\footnotesize Calea Grivi\c tei 21, P.O. Box 1-462, Bucharest, Romania\\[0.1cm]
\footnotesize E-mails: kohlenbach,leustean@mathematik.tu-darmstadt.de
}
\date{}
\maketitle

\begin{abstract}
This paper provides a fixed point theorem for asymptotically nonexpansive 
mappings in uniformly convex hyperbolic spaces as well as 
new effective results on the Krasnoselski-Mann iterations of such 
mappings. The latter were found using methods from logic and the paper 
continues a case study in the general program of extracting effective 
data from prima-facie ineffective proofs in the fixed point theory of 
such mappings.
\end{abstract}

\section{Introduction}

This paper provides a fixed point theorem for asymptotically nonexpansive 
mappings in uniformly convex hyperbolic spaces 
(Theorem \ref{FPP-ass-ne-uchyp}) as well as 
new effective results on the Krasnoselski-Mann iterations of such 
mappings (Theorem \ref{Herbrand-ass-ne}). 
The fixed point theorem generalizes corresponding theorems 
for uniformly convex normed spaces (\cite{Goebel+Kirk-72}) and CAT(0)-spaces 
(\cite{Kirk-2004}) while the effective bounds on the Krasnoselski-Mann 
iterations generalize results from \cite{K+Lambov} for the normed case  
which were obtained using techniques from mathematical logic or, more 
specifically, 
a proof theoretic method called 
(monotone) functional interpretation (see 
\cite{Kohlenbach(A),Gerhardy/Kohlenbach3}). In this respect the current
paper continues a case study in the general program of `proof mining' 
which is concerned with the extraction of effective uniform bounds from 
(prima-facie) ineffective proofs (see the discussion in 
section \ref{section-logic} and \cite{Kohlenbach(06a)} for a survey as 
well as \cite{Kohlenbach(book)}). 
Monotone functional interpretation 
systematically transforms any statement in a given proof into a new 
constructive version for which explicit bounds are provided. In the case of 
convergence statements (which this paper is about) 
this coincides with what recently has been 
advocated under the name `metastability' or `finite convergence' in 
an essay posted by T. Tao (\cite{Tao(07)}, see also 
\cite{Tao(07a)}). Thus the paper can also be seen as 
an instance of `hard analysis' as proposed by Tao. 
\\[2mm] 
Since the fundamental paper \cite{Goebel+Kirk-72}, the class of asymptotically 
nonexpansive mappings has been much studied in fixed point theory. Let 
$(X,d)$ be a metric space. A function $T:X\to X$ is called asymptotically 
nonexpansive if for some sequence $(k_n)$ in $[0,\infty )$ with 
$\lim\nolimits_{n\to\infty} k_n =0$ one has 
\[ d(T^n x,T^n y) \le (1+k_n)d(x,y), \ \ \forall n\in\N,\forall x,y\in X. \]
Asymptotically nonexpansive mappings have been studied mostly in the 
context of {\it uniformly convex} 
normed spaces (in fact for general normed spaces 
it is even open whether asymptotically nonexpansive selfmappings of 
bounded, closed, convex subsets have approximate fixed points, see 
\cite{Goebel(02)}). One typical result is the following theorem which 
is proved in \cite[Corollary 8]{K+Lambov} (as 
corollary of a quantitative result) but essentially is contained already  
in \cite{Rhoades(94),Schu(91),Schu(91B),Qihou(02)}):  
\begin{theorem} Let $(X,\|\cdot\|)$ be a uniformly convex normed space,  
$C\subseteq X$ a convex subset and $T:C\to C$ an asymptotically nonexpansive 
mapping with sequence $(k_n)$ in $[0,\infty )$ satisfying 
$\sum_{i=0}^{\infty} k_i <\infty.$ Let $(\lambda_n)$ be a sequence in $[a,b]$ 
for $0<a<b<1$ and define the Krasnoselski-Mann iteration of $T$ starting 
from $x\in X$ by 
\[ x_0:=x, \ x_{n+1} :=(1-\lambda_n)x_n+\lambda_n T^n(x_n). \] 
If $T$ has a fixed point, then $d(x_n,T(x_n))\stackrel{n\to\infty}{\to} 
0.$
\end{theorem}
While there does not seem to exist a computable rate of convergence 
in this case (in \cite{K+Lambov} it is shown that the proof 
even holds for asymptotically weakly-quasi nonexpansive functions
for which one can show that no uniform effective rate does exist), general 
logical metatheorems from \cite{Kohlenbach(metapaper),Gerhardy/Kohlenbach3} 
guarantee (see also section \ref{section-logic} below) 
effective uniform bound on the so-called no-counterexample 
interpretation of the convergence, or -- to use Tao's \cite{Tao(07),Tao(07a)} 
terminology -- 
on the metastability of $(\| x_n-T(x_n) \|),$ i.e. on 
\[ (*) \ \forall \varepsilon >0\,\forall g:\N\to\N \,\exists N\in\N 
\,\forall m\in [N,N+g(N)] \,( \| x_m-T(x_m)\| <\varepsilon), \]
which (ineffectively) is equivalent to the regular formulation of 
convergence towards $0.$ Here $[n,n+m]:=\{ n,n+1,n+2,\ldots,n+m\}.$ \\ 
The proof analyzed in \cite{K+Lambov} uses a lemma from 
\cite{Qihou(01A)}:
\begin{lemma}[\cite{Qihou(01A)}] \label{Qihou}
Let $(a_n),(b_n),(c_n)$ be sequences in $\R_+$ such that $\sum b_n$ and 
$\sum c_n$ are bounded and  
\[ \forall n\in\N (a_{n+1} \le (1+b_n)a_n +c_n). \] Then $(a_n)$ is 
convergent.
\end{lemma} 
The results in \cite{K+Lambov} were obtained by transforming a proof of 
$\| x_n-T(x_n)\| \to 0$ based on lemma \ref{Qihou} into a proof of $(*)$ 
together with an explicit effective bound for $(*)$ using a 
corresponding effective 
bound for the `metastability'-version of lemma \ref{Qihou} 
(see also proposition \ref{prop-quant-qihou-2-seq-v2} below) which 
constitutes a generalization of Tao's finite convergence principle from 
\cite{Tao(07)}. 
\\[2mm] 
In this paper we take the proofs from \cite{K+Lambov} as our 
point of departure and generalize the results to uniformly convex 
hyperbolic spaces (see the next section). This, in particular, covers the 
important class of CAT(0)-spaces (in the sense of Gromov) and, a-fortiorily, 
$\R$-trees in the sense of Tits. For CAT(0)-spaces we get a quadratic 
bound on the approximate fixed point property of $(x_n)$ (see 
corollary \ref{bounded-C-CAT0}).

\section{Hyperbolic spaces - definitions and properties}\label{hyp-spaces}

\noindent One can find in the literature  different notions of 'hyperbolic space' \cite{Kirk-1982,Goebel+Kirk-1983,Goebel+Reich-book,Reich+Shafrir-1990}. We work in the setting of hyperbolic spaces as introduced by the first author \cite{Kohlenbach(metapaper)}, which are slightly more restrictive than the spaces of hyperbolic type in the sense of Goebel/Kirk \cite{Goebel+Kirk-1983}, but more general than the hyperbolic spaces in the sense of Reich/Shafrir \cite{Reich+Shafrir-1990}.

A {\em  hyperbolic space}  $(X,d,W)$ is a metric space $(X,d)$ together with a convexity mapping $W:X\times X\times [0,1]\to X$ satisfying 
\begin{eqnarray*}
(W1) & d(z,W(x,y,\lambda))\le (1-\lambda)d(z,x)+\lambda d(z,y),\\
(W2) & d(W(x,y,\lambda),W(x,y,\tilde{\lambda}))=|\lambda-\tilde{\lambda}|\cdot 
d(x,y),\\
(W3) & W(x,y,\lambda)=W(y,x,1-\lambda),\\
(W4) & \,\,\,d(W(x,z,\lambda),W(y,w,\lambda)) \le (1-\lambda)d(x,y)+\lambda
d(z,w).
\end {eqnarray*}

The convexity mapping $W$ was first considered by Takahashi in \cite{Takahashi-70}, where a triple $(X,d,W)$ satisfying $(W1)$ is called a convex metric space.

The class of hyperbolic spaces includes normed spaces and convex subsets thereof, the Hilbert ball \cite{Goebel+Reich-book} as well as CAT(0)-spaces in the sense of Gromov (see \cite{Bridson+Haefliger-book} for a detailed treatment).

\noindent If $x,y\in X$ and $\lambda\in[0,1]$ then we use the notation $(1-\lambda)x\oplus \lambda y$ for $W(x,y,\lambda)$.  It is easy to see that for any $x,y\in X$ and any $\lambda\in[0,1]$, 
\begin{equation}
d(x,\lambdaxy)=\lambda d(x,y),\text{~and~} d(y,\lambdaxy)=(1-\lambda)d(x,y).\label{prop-xylambda}
\end{equation}

\noindent  We shall denote by $[x,y]$ the set $\{(1-\lambda)x\oplus \lambda y:\lambda\in[0,1]\}$. A nonempty subset $C\subseteq X$ is {\em convex} if $[x,y]\in C$ for all $x,y\in C$. 

For any $x\in X, r>0$, the open (closed) ball with center $x$ and radius $r$ is denoted with $U(x,r)$ (respectively $\ol{U}(x,r)$). It is easy to see that open and closed balls are convex. Moreover, using (W4), we get that the closure of a convex subset of a hyperbolic spaces is again convex.

One of the most important classes of Banach spaces are the uniformly convex ones, introduced by Clarkson in the 30's \cite{Clarkson-1936}. 
Following \cite[p. 105]{Goebel+Reich-book}, we can define uniform convexity for hyperbolic spaces too.

A hyperbolic space $(X,d,W)$ is  {\em uniformly convex} \cite{L-07-JMAA} if for 
any $r>0$ and any $\varepsilon\in(0,2]$ there exists $\theta\in(0,1]$ such that 
for all $a,x,y\in X$,
\begin{eqnarray}
\left.\begin{array}{l}
d(x,a)\le r\\
d(y,a)\le r\\
d(x,y)\ge\varepsilon r
\end{array}
\right\}
& \quad \Rightarrow & \quad d\left(\frac12x\oplus\frac12y,a\right)\le (1-\theta)r. \label{uc-def}
\end{eqnarray}
A mapping $\eta:(0,\infty)\times(0,2]\rightarrow (0,1]$ providing such a
$\theta:=\eta(r,\varepsilon)$ for given $r>0$ and $\varepsilon\in(0,2]$ is called a {\em modulus of uniform convexity}. 

In the sequel, $(X,d,W)$ is a uniformly convex space and $\eta$ is a modulus of uniform convexity.

\blem\label{eta-prop-1}
Let $r>0,\varepsilon\in(0,2]$ and $a,x,y\in X$ be such that $d(x,a)\le r,d(y,a)\le r ,d(x,y)\ge\eps r$. Then for any $\lambda\in[0,1]$,
\be
\item \label{uc-ineq-Groetsch} $\ds d(\lambdaxy,a)\le  (1-2\lambda(1-\lambda)\eta(r,\varepsilon))r$; 
\item for any $\psi\in (0,2]$ such that $\psi\le\eps$, 
\[\ds d(\lambdaxy,a)\le  (1-2\lambda(1-\lambda)\eta(r,\psi))r\,;\]
\item  \label{eta-s-geq-r} for any $s\geq r$, 
\[d(\lambdaxy,a) \le \left(1-2\lambda(1-\lambda)\eta\left(s,\eps\frac{r}{s}\right)\right)s\,.\]
\ee
\elem
\begin{proof}
\be
\item See \cite[Lemma 7]{L-07-JMAA}.
\item Note that $d(x,y)\ge\eps r\ge\psi r$ and apply \ref{uc-ineq-Groetsch}.\item Since  $d(x,a),d(y,a)\le r\le s$, $\ds d(x,y)\ge\eps r=\left(\eps\frac{r}{s}\right)s$ and $\ds 0<\eps\frac{r}{s}\le \eps\leq 2$, the conclusion follows again by an application of \ref{uc-ineq-Groetsch}.
\ee
\end{proof}

\noindent We say that $\eta$ is {\em monotone} if it decreases with $r$ (for a fixed $\epsilon$). It turns out that CAT(0)-spaces  are  uniformly convex hyperbolic spaces having a monotone modulus of uniform convexity, quadratic in $\eps$:  $\ds\eta(r,\varepsilon)=\varepsilon^2/8$.  We refer to \cite{L-07-JMAA} for details.

The following proposition is one of the main ingredients in the proof of Theorem \ref{FPP-ass-ne-uchyp}. Its proof is similar to the one of the corresponding result for uniformly convex Banach spaces (see, for example, \cite[Theorem 2.1]{Goebel+Reich-book}).
\bprop\label{CIP}
Let $(X,d,W)$ be a complete uniformly convex hyperbolic space with a monotone modulus of uniform convexity $\eta$. \\
The intersection of any decreasing sequence of nonempty bounded closed convex subsets of $X$ is nonempty.
\eprop
\begin{proof}
Let $(C_n)_{n\geq 1}$ be a decreasing sequence of nonempty bounded closed convex subsets of $X$ and let $x\in X$ be arbitrary. If $x\in C_n$ for all $n\in\N$, then $\ds\bigcap_{n\geq 1} C_n\ne\emptyset$. Assume that there exists $N\in\N$ such that $x\nin C_N$, so that $d(x,C_N)>0$, since $C_N$ is closed. If $r_n:=d(x,C_n)$, then $(r_n)$ is an increasing sequence of nonnegative reals, bounded from above by $d(x,a)+diam(C_1)$, where $a\in C_1$. It follows that $r:=\lim r_n =\sup r_n \geq r_N>0$.

Define $\ds D_n:=C_n\cap \ol{U}\left(x,r+\frac1n\right)$.  Then it is easy to see that $(D_n)$ is a decreasing sequence of nonempty closed subsets of $X$. Let $d_n:=diam(D_n)$ and $0\leq d:=\lim d_n=\inf d_n$. 

Assume that $d>0$. Let $K\in \N$ be such that $\ds\frac1K\leq \frac{d}2$. For any $n\geq K$, there exist $x_n,y_n\in D_n$ such that $\ds d(x_n,y_n)\geq d_n-\frac1n\geq d -\frac1n\geq \frac{d}2$. 

Since $\ds d(x_n,x), d(y_n,x)\le r+\frac1n,\,\ds d(x_n,y_n)\geq \frac{d}2\geq \left(r+\frac1n\right)\cdot \frac{d}{2(r+1)}$ and $\ds\frac{d}{2(r+1)}\leq 1$, we get that for all $n\geq K$,
\bua
r_n&\leq& d\left(\frac12x_n\oplus\frac12y_n,x\right)\le \left(1-\eta\left(r+\frac1n ,\frac{d}{2(r+1)}\right)\right)\cdot\left(r+\frac1n\right), \\
&& \text{since~~} X \text{~~is uniformly convex}\\
&\leq& \left(1-\eta\left(r+1,\frac{d}{2(r+1)}\right)\right)\cdot\left(r+\frac1n\right),\\
&& \text{since~~} r+\frac1n\leq r+1 \text{~~and~~} \eta \text{~~is monotone.}
\eua
Thus, by letting $n\to\infty$, $\ds r\leq \left(1-\eta\left(r+1,\frac{d}{2(r+1)}\right)\right)\cdot r< r$, that is a contradiction.  \\
It follows that we must have $d=0$. This and the completeness of $X$ imply that $\ds\bigcap_{n\geq 1} D_n\ne\emptyset$, hence $\ds\bigcap_{n\geq 1} C_n\ne\emptyset$.
\end{proof}

\section{Main results}\label{main-results}

The notion of nonexpansive mapping can be introduced in the very general setting of metric spaces. Thus, if $(X,d)$ is a metric space, and $C\subseteq X$ a nonempty subset, than a mapping $T:C\to C$ is called {\em nonexpansive} if for all $x,y\in C$,
\[d(Tx, Ty)\le d(x,y).\]

Asymptotically nonexpansive mappings were introduced  by Goebel and Kirk \cite{Goebel+Kirk-72} as a generalization of the nonexpansive ones. A function $T:C\to C$ is said to be {\em asymptotically
nonexpansive with sequence}  $(k_n)_{n\geq 0}$ in $[0,\infty)$ if 
$\lim\limits_{n\to\infty} k_n =0$ and 
\[  d(T^n x,T^ny) \leq (1+k_n)d(x,y), \ \hfill\forall n\in\N,
\forall x,y\in C.
 \] 
$Fix(T)$ denotes the set of fixed points of $T$ and for any $\eps>0$, $Fix_\eps(T)$ denotes the set of $\eps$-fixed points, that is points $x\in C$ such that $d(x,Tx)<\eps$.

We say that $C$ has the {\em fixed point property} (FPP)   for asymptotically nonexpansive mappings if $Fix(T)\ne\emptyset$ for any asymptotically nonexpansive mapping $T:C\to C$. Moreover, $C$ has the {\em approximate fixed point property} (AFPP)   for asymptotically nonexpansive mappings if $Fix_\eps(T)\ne\emptyset$ for any asymptotically nonexpansive mapping $T:C\to C$ and any $\eps>0$.

Goebel and Kirk proved the following generalization of the famous Browder-Goehde-Kirk fixed point  theorem for nonexpansive mappings.

\bthm\cite[Theorem 1]{Goebel+Kirk-72} \label{FPP-ucBanach}\\
Nonempty closed convex and bounded subsets of uniformly convex Banach spaces have the FPP for asymptotically nonexpansive mappings.
\ethm

\noindent In 2004, Kirk obtained a similar result for CAT(0)-spaces.

\bthm\cite[Theorem 28]{Kirk-2004} \label{FPP-CAT0}\\
Nonempty closed convex and bounded subsets of complete CAT(0)-spaces have the FPP for asymptotically nonexpansive mappings.
\ethm

\noindent Kirk proved Theorem \ref{FPP-CAT0} using nonstandard methods, inspired by Khamsi's proof that bounded hyperconvex metric spaces have the AFPP for asymptotically nonexpansive mappings \cite{Khamsi-2003}.

The first main result of this paper is a generalization of Theorem \ref{FPP-ucBanach} to uniformly convex hyperbolic spaces with monotone modulus of uniform convexity.

\bthm\label{FPP-ass-ne-uchyp}
Let $(X,d,W)$ be a complete uniformly convex hyperbolic space having a monotone modulus of uniform convexity. Then any nonempty closed convex and bounded subset of $X$ has the FPP for asymptotically nonexpansive mappings.
\ethm

\noindent Our proof follows closely Goebel and Kirk's proof of Theorem \ref{FPP-ucBanach} and we present the details in Section \ref{FPP-uchyp-proof}. As a consequence, we obtain also an elementary proof of Theorem \ref{FPP-CAT0}. 

In fact, as it was already pointed out for uniformly convex normed spaces in \cite{K+Lambov}, the proof of the FPP can be transformed into an elementary proof of the AFPP,  which does not need the completeness of $X$ or the closedness of $C$.

\bprop\label{AFPP-ass-ne-uchyp}
Let $(X,d,W)$ be a uniformly convex hyperbolic space having a monotone modulus of uniform convexity. Then any nonempty  convex and bounded subset of $X$ has the AFPP for asymptotically nonexpansive mappings.
\eprop
\begin{proof}
The proof of \cite[Lemma 21]{K+Lambov} generalizes easily to our setting.
\end{proof}

The main part of the paper will be devoted to getting a quantitative version  of an asymptotic regularity theorem of the Krasnoselskii-Mann iterations of asymptotically nonexpansive mappings.

Let $(X,d,W)$ be a hyperbolic space, $C\se X$ a nonempty convex subset of $X$ and $T:C\to C$ an asymptotically nonexpansive mapping. 

For asymptotically nonexpansive mappings, the {\em Krasnoselski-Mann iteration} starting from $x\in C$ is defined by:   
\begin{equation}
x_0:=x, \quad x_{n+1}:=(1-\lambda_n)x_n \oplus\lambda_n T^nx_n, \label{KM-lambda-n-def-hyp}\end{equation}
where $(\lambda_n)$ is a sequence in $[0,1]$.

Following \cite{Borwein+Reich+Shafrir-92}, we say that $T$ is {\em $\lambda_n$-asymptotically regular} if  for all $x\in C$, 
\[\lim_{n\to\infty}d(x_n,Tx_n)=0.\]

The second main result of the paper is the following theorem, generalizing to uniformly convex hyperbolic spaces a similar result obtained for uniformly convex normed spaces by the first author and Lambov \cite{K+Lambov}. 

\begin{theorem} \label{Herbrand-ass-ne}
Let $(X,d,W)$ be a uniformly convex hyperbolic space with a monotone modulus of uniform convexity $\eta$, $C$ be a nonempty convex subset of $X$ and $T:C\to C$ be asymptotically nonexpansive with sequence $(k_n)$. \\
Assume that $K\geq 0$ is such that $\ds\sum_{n=0}^\infty k_n\le K$ and that $L\in\N, L\geq 2$  is such that $\ds\frac{1}{L}\leq \lambda_n \leq 1-\frac{1}{L}$ for all $n\in\N$.

Let $x\in C$ and $b>0$ be such that for any $\delta>0$ there is $p\in C$ with
\beq
d(x,p)\leq b\si d(Tp,p)\leq \delta.\label{hyp-app-b-x}
\eeq
Then for all $\eps\in (0,1]$ and for all $g:\N\to\N$,
\beq
\exists N\le \Phi(K, L, b, \eta, \eps, g)
\forall m\in [N, N+g(N)] \left (d(x_m,Tx_m) < \eps \right), \label{Herbrand-ass-ne-conclusion}
\eeq 
where 
\[\ba{l}
\Phi(K, L, b, \eta, \eps, g):=\ds h^M(0), \quad h(n) := g(n+1)+n+2,\\[0.1cm]
M :=\left\lceil\ds\frac {3\left(5KD+D+\frac{11}{2}\right)}{\theta}\right\rceil, \quad  
\ds D := e^K\left(b+ 2\right),\\
\theta := \ds\frac{\eps}{L^2f(K)}\cdot\eta\left((1+K)D+1,\frac{\eps}{f(K)((1+K)D+1)}\right),\quad \\
f(K) := 2(1+(1+K)^2(2+K)).
\ea\]
Moreover, $N=h^i(0)+1$ for some $i<M$.
\end{theorem}

\noindent We shall give the proof of the above theorem in the last section of our paper. As we shall explain in detail in Section \ref{section-logic}, the extractability of the bound $\Phi$ is guaranteed by a general logical metatheorem. Moreover, this theorem allows us to conclude that $\ds \lim d(x_n,Tx_n)=0$, assuming the existence of approximate fixed points in some neighborhood of the starting point $x\in C$ (see  the discussion on the Herbrand normal form in Section \ref{section-logic}). 

\begin{fact}\label{theorem-weaken-hyp}
By an inspection of its proof, it is easy to see that the above theorem remains true if we weaken the hypotheses on $(k_n)$ and $(\lambda_n)$. In fact, it is enough to require that $\ds\sum_{n=0}^\Phi k_n\le K$ and  $\ds\frac{1}{L}\leq \lambda_n \leq 1-\frac{1}{L}$ for all $n\leq \Phi$. Note that once the hypotheses are weakened one must move these hypotheses under the scope of the quantification over $\eps$ and $g$ since $\Phi$ depends on these.
\end{fact}

\begin{fact}\label{theorem-tilde-eta}
Assume, moreover, that  $\eta(r,\varepsilon)$ can be written as $\eta(r,\varepsilon)=\varepsilon\cdot\tilde{\eta}(r,\varepsilon)$ such that $\tilde{\eta}$ increases with $\varepsilon$ (for a fixed $r$). Then we can replace $\eta$ with $\tilde{\eta}$ in the bound $\Phi(K, L, b, \eta, \eps, g)$.
\end{fact}
\begin{proof}
Define 
\[\theta := \ds\frac{\eps}{L^2f(K)}\cdot\tilde{\eta}\left((1+K)D+1,\frac{\eps}{f(K)((1+K)D+1)}\right)\]
and follow the proof of the theorem using Lemma \ref{main-technical-lemma}, (\ref{eta-tilde}) instead of Lemma \ref{main-technical-lemma}, (\ref{eta-general}).
\end{proof}

We give now some further corollaries.

\begin{theorem} \label{AsReg-as-ne}
Assume $(X,d,W),\eta, C, T:C\to C, (k_n), K, (\lambda_n), L$ are  as in the hypotheses of Theorem \ref{Herbrand-ass-ne}.

Let  $x\in C$ and $b>0$ be such that for any $\delta>0$ there is $p\in C$ with
\beq
d(x,p)\leq b\si d(Tp,p)\leq \delta.
\eeq
Then $\ds \lim d(x_n,Tx_n)=0$ and, moreover,
\beq
\forall\eps\in (0,1]\exists N\le \Phi(K, L, b, \eta, \eps)
\left (d(x_N,Tx_N) \leq \eps \right),
\eeq 
where $\Phi(K, L, b, \eta, \eps):=\ds 2M$ and $M, D, \theta, f(K)$ are as in Theorem \ref{Herbrand-ass-ne}. 
\end{theorem}
\begin{proof}
Take $g(n)\equiv 0$ in Theorem \ref{Herbrand-ass-ne}.
\end{proof}

\bcor (see also Theorem \ref{original-theorem})\\
Assume $(X,d,W),\eta, C, T:C\to C, (k_n), K, (\lambda_n), L$ are  as in the hypotheses of Theorem \ref{Herbrand-ass-ne}.

If $Fix(T)\ne\emptyset$, then $T$ is $\lambda_n$-asymptotic regular.
\ecor
\begin{proof}
Let $\tilde{p}$ be a fixed point of $T$. For any $x\in C$,  (\ref{hyp-app-b-x}) is satisfied with $b:=d(x,\tilde{p})$ and $p:=\tilde{p}$. 
\end{proof}

\begin{corollary}\label{bounded-C}
Let $(X,d,W),\eta, C, T:C\to C, (k_n), K, (\lambda_n), L$ be  as in the hypotheses of Theorem \ref{Herbrand-ass-ne}. Assume moreover that $C$ is bounded with finite diameter $d_C$.

Then $T$ is $\lambda_n$-asymptotic regular, and  the following holds for all $x\in C$:
\beq
\forall\eps\in (0,1]\exists N\le \Phi(K, L, d_C, \eta, \eps)
\left (d(x_N,Tx_N) < \eps \right),
\eeq 
where $\Phi(K, L, d_C, \eta, \eps)$ is defined as in Theorem \ref{AsReg-as-ne} by replacing $b$ with $d_C$.
\end{corollary}
\begin{proof}
If $C$ is bounded, then $C$ has the AFPP for asymptotically nonexpansive mappings by Proposition \ref{AFPP-ass-ne-uchyp}, so the condition (\ref{hyp-app-b-x}) holds for all $x\in C$ with $d_C$ instead of $b$. Hence, we can  conclude that $\ds \lim d(x_n,Tx_n)=0$ for all $x\in C$.
\end{proof}

\noindent Thus, for bounded $C$, we get asymptotic regularity and an explicit approximate fixed point bound $\Phi(K, L, d_C, \eta, \eps)$, which depends only on the error $\varepsilon$, on the modulus of uniform convexity $\eta$, on the diameter $d_C$ of $C$, on $(\lambda_n)$ via $L$ and on $(k_n)$ via $K$, but not  on the nonexpansive mapping $T$, the starting point $x\in C$ of the iteration or other data related with $C$ and $X$.

As we have pointed out in Section \ref{hyp-spaces}, CAT(0)-spaces are uniformly convex hyperbolic spaces with a 'nice' monotone modulus of uniform convexity $\ds \eta(r,\varepsilon):=\frac{\varepsilon^2}{8}$. Hence, as an immediate consequence of Corollary \ref{bounded-C} and Remark \ref{theorem-tilde-eta} we get the following result. 

\begin{corollary}\label{bounded-C-CAT0}
Let $X$ be a CAT(0)-space, $C$ be a nonempty convex bounded subset of $X$ with diameter $d_C$ and $T:C\to C$ be asymptotically nonexpansive with sequence $(k_n)$. \\
Assume that $K\geq 0$ is such that $\ds\sum_{n=0}^\infty k_n\le K$ and that $L\in\N, L\geq 2$  is such that $\ds\frac{1}{L}\leq \lambda_n \leq 1-\frac{1}{L}$ for all $n\in\N$.

Then $T$ is $\lambda_n$-asymptotic regular, and  the following holds for all $x\in C$:
\beq
\forall\eps\in (0,1]\exists N\le \Phi(K, L, d_C,\eps)
\left (d(x_N,Tx_N) < \eps \right),
\eeq 
where 
\[\ba{l}
\Phi(K, L, d_C, \eps):=2M, \\
M :=\left\lceil\ds\frac{1}{\eps^2}\cdot 24L^2\left(5KD+D+\frac{11}{2}\right)(f(K))^3((1+K)D+1)^2\right\rceil, \\
  D := \ds e^K\left(d_C+ 2\right), \quad f(K) := 2(1+(1+K)^2(2+K)).
\ea\]
\end{corollary}

\noindent Hence, in the case of convex bounded subsets of CAT(0)-spaces, we get a quadratic (in $1/\varepsilon$) approximate fixed point bound.  We recall that  for nonexpansive mappings, a quadratic rate of asymptotic regularity for the Krasnoselski-Mann iterations was obtained by the second author  \cite{L-07-JMAA}.

\section{Proof of Theorem \ref{FPP-ass-ne-uchyp}}\label{FPP-uchyp-proof}

In this section, we give the proof of Theorem \ref{FPP-ass-ne-uchyp}. As we have already pointed out,  we generalize to our setting Goebel and Kirk's proof for uniformly convex Banach spaces. \\

\noindent{\em Proof of Theorem \ref{FPP-ass-ne-uchyp}}

\noindent For any $y\in C$, let us consider
\bua
A_y:=\left\{ a\in \R_+ \mid  \text{there exist~}  x \in C, k\in\N
\text{~such that~ } d(T^iy, x) \le a  \text{~for all~}i \geq k\right \}.
\eua
If $d(C)$ is the diameter of $C$, then  $d(C)\in A_y$,  hence $A_y$ is nonempty. Let $\alpha_y:=\inf A_y$. For any $\theta>0$ there exists $a_\theta\in A_y$ such that $\ds a_\theta <\alpha_y+\theta$, so
\beq
\exists x \in C \exists k \in \N \forall i\geq  k\left(d(T^iy,x) \leq a_\theta< \alpha_y+ \theta\right).\label{alpha-y-prop}
\eeq

\noindent Obviously, $\alpha_y\geq 0$. We distinguish two cases:

\noindent {\bf Case 1.} $\alpha_y=0$. 

\noindent Let $\eps>0$. Applying  (\ref{alpha-y-prop}) with $\ds\theta:=\frac\eps2$, we get the  existence of $x\in C$ and $k\in\N$ such that  for all $m,n\geq k$
\beq
d(T^my,T^ny)\leq d(T^my,x)+d(T^ny,x)<\frac{\eps}{2}+\frac{\eps}{2}=\eps,
\eeq
so the sequence $(T^ny)_{n\geq 1}$ is Cauchy, hence convergent to some $z\in C$. It is easy to see that $z$ is a fixed point of $T$.  

\noindent {\bf Case 2.} $\alpha_y>0$. 

For any $n\geq 1$, let us define 
\bea
C_n:= \bigcup_{k\geq 1}\bigcap_{i\geq k}\ol{U}\left(T^iy,\alpha_y+\frac1n\right), & D_n:=\ol{C_n}\cap C.
\eea
By (\ref{alpha-y-prop}) with $\ds\theta:=\frac1n$, there exist $x\in C, k\geq 1$ such that $\ds x\in \bigcap_{i\geq k}\ol{U}\left(T^iy,\alpha_y+\frac1n\right)$, hence $D_n$ is nonempty. Moreover, $(D_n)_{n\geq 1}$ is a decreasing sequence of nonempty bounded closed convex subsets of $X$, hence we can apply Proposition \ref{CIP} to get that 
 \[D:=\ds\bigcap_{n\geq 1}D_n\ne\emptyset.\]

\noindent {\bf Claim:} For any $x\in D$ and $\theta>0$ there exists $K\in\N$ such that for all $i\geq K$,
\beq d(T^iy,x)\leq \alpha_y+\theta.\label{claim}
\eeq
{\bf Proof of claim:} Let $x\in D,\theta>0$ and $N\in\N$ be such that $\ds\frac2N\leq\theta$. Since $x\in D$, we have that $x\in\ol{C_N}$, so  there exists a sequence $(x^N_n)_{n\geq 1}$ in $C_N$ such that $\lim x^N_n=x$. Let $P\geq 1$ be such that $\ds d(x,x^N_n)\leq \frac1N$ for all $n\geq P$ and $K\geq 1$ such that $x^N_P\in \bigcap_{i\geq K}\ol{U}\left(T^iy,\alpha_y+\frac1N\right)$.

It follows that for all $i\geq K$,
\bua
d(T^iy,x)\leq d(T^iy,x^N_P)+d(x^N_P,x)\leq  \alpha_y+\frac1N+\frac1N=
\alpha_y+\frac2N\leq \alpha_y+\theta.
\eua
Thus, the claim is proved.

In the sequel, we shall prove that any point of $D$ is a fixed point of $T$. Let $x\in D$ and assume by contradiction that $Tx\ne x$. Then 
$(T^nx)$ does not converge to $x$, so there exists $\eps>0$ such that
\beq 
\forall k \in \N\exists n\geq k(d(T^nx,x)\geq \varepsilon/2).\label {fixe_xdelta_not_fixe}
\eeq
We can of course assume that $\eps\in(0,4]$. Then $\ds\frac{\varepsilon}{2(\alpha_y + 1)}\in (0,2]$ and there exists $\theta_y\in(0,1]$ such that
\beq
1-\eta \left(\alpha_y+1,\frac{\varepsilon}{2(\alpha_y + 1)}\right) \leq \frac{\alpha_y - \theta_y}{\alpha_y+\theta_y}.\label{def-delta-y}
\eeq 
Since $\ds \lim (1+k_n)\left(\alpha_y+\frac{\theta_y}{2}\right)=\alpha_y+\frac{\theta_y}{2}<\alpha_y+\theta_y$, there exists $N_0\in\N$ such that
\beq
\forall n\geq N_0\left((1+k_n)\left(\alpha_y+\frac{\theta_y}{2}\right)<\alpha_y+\theta_y\right).\label{hyp-k-n}
\eeq
Applying (\ref{claim}) with $\ds\theta:=\frac{\theta_y}2$, there exists $K\in\N$ such that
\beq
\forall i\geq K\left(d(T^iy,x)\leq \alpha_y+\frac{\theta_y}{2}\right). \label{def-z-K}
\eeq
Applying (\ref{fixe_xdelta_not_fixe}) with $k:=N_0$, we get $N\geq N_0$ such that
\beq
d(T^{N}x,x)\geq \varepsilon/2.\label{hyp1-N}
\eeq
Let now $m\in \N$ be such that $m\geq N+K$.  Then
\bua
d(T^Nx, T^my) &=& d(T^Nx, T^N(T^{m-N}y))\leq (1+k_N)d(x, T^{m-N}y)\\
&<& (1+k_N)\left(\alpha_y+\frac{\theta_y}{2}\right),\quad \text{~by~} (\ref{def-z-K})\\
&< & \alpha_y+\theta_y, \quad \text{~by~} (\ref{hyp-k-n}).
\eua
Hence,
\bua
d(T^Nx, T^my) &< & \alpha_y+\theta_y, \\
d(x, T^my)&\leq & \alpha_y+\frac{\theta_y}{2}<\alpha_y+\theta_y, \quad \text{~by~} (\ref{def-z-K}),\\
d(x,T^Nx) &\geq &\frac{\varepsilon}{2}, \quad \text{~by~} (\ref{hyp1-N})\\
&=& (\alpha_y+\theta_y)\frac{\eps}{2(\alpha_y+\theta_y)}\geq (\alpha_y+\theta_y)\frac{\eps}{2(\alpha_y+1)}.
\eua
Applying now the fact that $X$ is uniformly convex, we get that
\bua
\rho\left(\frac12 x\oplus\frac12 T^Nx,T^my\right)&\leq &
\left(1-\eta\left(\alpha_y+\theta_y,\frac{\eps}{2(\alpha_y+1)}\right)\right)(\alpha_y+\theta_y).
\eua
Since  $\alpha_y+\theta_y\leq \alpha_y+1$ and  $\eta$ is monotone,
\bua
1-\eta\left(\alpha_y+\theta_y,\frac{\eps}{2(\alpha_y+1)}\right)&\leq&1- \eta\left(\alpha_y+1,\frac{\eps}{2(\alpha_y+1)}\right)
\leq \frac{\alpha_y - \theta_y}{\alpha_y+\theta_y}, \quad\text{by~~} (\ref{def-delta-y}).
\eua
Hence,
\bua
\rho\left(\frac12 x\oplus\frac12 T^Nx,T^my\right)&\leq & \frac{\alpha_y - \theta_y}{\alpha_y+\theta_y}\cdot(\alpha_y+\theta_y)=\alpha_y - \theta_y.
\eua

Thus, there exist $k:=N+K$ and $\ds z:=\frac12 x\oplus\frac12 T^Nx \in C$ such that for all $m\geq k$, $d(z,T^my)\leq\alpha_y - \theta_y$. This means that $\alpha_y - \theta_y\in A_y$. Since $\alpha_y - \theta_y<\alpha_y=\inf A_y$, we have got a contradiction.\\
It follows that $x$ is a fixed point of $T$.
\qed

\section{A general logical metatheorem} \label{section-logic}

\noindent One of the main results of this paper, Theorem \ref{Herbrand-ass-ne},  is a quantitative version of an asymptotic regularity theorem for asymptotically nonexpansive mappings of hyperbolic spaces. In this section we indicate how such a version can be obtained 
from a prima-facie ineffectively proven convergence result by means of 
a general logical metatheorem. Such metatheorems were developed first in 
\cite{Kohlenbach(metapaper)} and \cite{Gerhardy/Kohlenbach3} and guarantee 
for general classes of theorems and proofs the extractability of uniform 
effective bounds from given proofs (see 
\cite{Kohlenbach(book)} for a comprehensive treatment). 
The metatheorems apply to general classes of 
spaces such as metric, hyperbolic, normed, uniformly convex and inner 
product spaces (as well as their completions) and functions such as 
nonexpansive, Lipschitz, weakly quasi-nonexpansive or uniformly 
continuous functions among others. We state here only one particular 
corollary of such a metatheorem which covers the situation treated 
in this paper.

The formal system ${\cal A}^{\omega}[X,d,W]_{-b}$ results from the extension 
of a system ${\cal A}^{\omega}$ for analysis (going back to Spector 
\cite{Spector(62)}) obtained by axiomatizing an abstract hyperbolic space 
$(X,d,W).$ 
This is achieved by adding constants $d_X$ and $W_X$ representing 
$d,W$ to the system together with axioms expressing that $d_X$ is a 
pseudo-metric and $W_X$ satisfies the axioms (W1)-(W4) 
(the subscript `$-b$' refers to the fact that we do not assume 
$(X,d,W)$ to be bounded). Equality for 
objects in $X$ is defined as $x=_X y:\equiv d_X(x,y)=_{\R} 0$ so that 
we actually consider the {\bf metric} 
space induced by the pseudo-metric $d_X.$ 
The language of ${\cal A}^{\omega}[X,d,W]_{-b}$ is based on the language of  
functionals of all types over $\N,X$ together with appropriate induction and 
recursion axioms as well as the axiom schema of dependent choice for all 
types (which, in particular, implies countable choice and -- as a consequence 
of this -- full comprehension over natural numbers). So in particular full 
so-called 2nd order arithmetic is a subsystem of ${\cal A}^{\omega}.$ 
Precise definitions for all this can be found in 
\cite{Kohlenbach(metapaper),Gerhardy/Kohlenbach3}. 
To have quantifiers for 
functionals over $\N,X$ means that we can quantify not only over 
$\N$ (starting from $0$) and $X$ but also over 
functions $f:\N\to\N,$ $g:X\to X$ 
$h:\N\to X$ (i.e. sequences in $X$) and even over function(al)s taking 
such objects as arguments and so on. The types $\N,\N\to\N$ 
(and also $k$-ary number-theoretic functions), $X,X\to X,\N\to X$ are called 
{\em small types}. Treating general so-called Polish (i.e. complete separable) metric 
spaces $P$ as continuous images of the Baire space $\N^{\N}$ the type 
$\N\to\N$ also covers quantification over $P$ (for Polish spaces $P$ 
given in so-called standard representation).  
${\cal A}^{\omega}[X,d,W,\eta]_{-b}$ results from 
${\cal A}^{\omega}[X,d,W]_{-b}$ by adding a new constant $\eta:\N\times\N
\to\N$ together with axioms expressing that $\eta$ represents a modulus of 
uniform convexity of $(X,d,W)$ (see \cite{L-2006}).
$A_{\forall}$ (resp. $B_{\exists})$ is called a $\forall$-formula 
(resp. an $\exists$-formula) if it has the form $\forall \underline{a}\,
A_{qf}(\underline{a})$ (resp. $\exists \underline{a}\,B_{qf}(\underline{a})$) 
where $\underline{a}$ is a tuple of variables of small types and $A_{qf}$ 
($B_{qf}$) is a quantifier-free formula. \\[2mm]     
Let $A_{\forall}(x,y,z,T,u)$ and $B_{\exists}(x,y,z,T,v)$ be 
$\forall$- resp. $\exists$-formulas  
which only contain the shown 
variables as free variables. In the following we abbreviate 
$A_{\forall}(x,y,z,T,u)$ and $B_{\exists}(x,y,z,T,v)$ by $A_{\forall}$ and 
$B_{\exists}.$ For $T:X\to X,$ $x\in X$ and $b\in\N$, the 
formula $Fix_{\varepsilon}(T,x,b)\not=\emptyset$ expresses that 
$T$ has an $\varepsilon$-fixed point $p$ in the $b$-ball around 
$x,$ i.e. $d(x,p)\le b$ and $d(p,T(p)) <\varepsilon.$
\begin{theorem}
\label{meta-corollary}
Let $(\lambda_n)$ be some standard enumeration of $\Q^*_+.$ 
  \begin{enumerate}
  \item (\cite{Gerhardy/Kohlenbach3}, Corollary 4.26) \\  
 Let $P$ (resp. $K$) be a ${\cal A}^{\omega}$-definable Polish space
  (resp. compact metric space).
Assume one can prove in ${\cal A}^{\omega}[X,d,W]_{-b}$ a sentence:
    \[ \ba{l} 
    \forall x\in P\ \forall y\in K\ \forall n\in\N \ \forall z\in X \ 
\forall T:X \to X \\ \hspace*{1cm}
(\mbox{$T\ \lambda_n$-Lipschitz  } \wedge Fix(T)\not=\emptyset\wedge 
    \forall u\in\N \ A_{\forall} \rightarrow \exists
    v\in\N \ B_{\exists} ).\ea
    \]
    Then from the proof one can extract a computable\footnote{Here we 
refer to the usual oracle version (`type-2') of computability when dealing 
with arguments in $\N^{\N}.$}  functional $\Phi: 
    \N^{\N}\times\N\times\N
    \to \N$ s.t. for all representatives $r_x \in\N^{\N}$ of $x\in P$ and 
    all $n,b\in\N$
    \[ \begin{array}{l}
    \forall y\in K\ \forall z\in X \ \forall T:X \to X \ (T \ \mbox
{$\lambda_n$-Lipschitz} \ \wedge\forall\varepsilon >0 
Fix_{\varepsilon}(T,x,b)\not=\emptyset \\ 
\hspace*{1cm}\wedge 
    \, d_X(z,T(z)) 
    \le_{\R} b \wedge   
    \forall u \leq \Phi(r_x,n,b) \ A_{\forall} 
    \rightarrow \exists v \leq \Phi(r_x,n,b) \ 
    B_{\exists})
    \end{array}
    \]
    holds in all (nonempty) hyperbolic spaces $(X,d,W)$. 
\item (\cite{L-2006}) 
If the premise of this rule is proved in  
${\cal A}^{\omega}[X,d,W,\eta]_{-b},$ then the conclusion 
holds in all (nonempty) uniformly convex hyperbolic spaces $(X,d,W)$ 
provided that $\eta$ is interpreted by a modulus of uniform 
convexity of $(X,d,W).$ The bound $\Phi$ then additionally depends on 
$\eta.$
\end{enumerate} 
For the special cases where $P$ is $\N$ with the discrete metric resp. 
$\N^{\N}$ with the product metric (Baire space), we can treat the elements 
of $P$ directly without any representation, i.e. $r_x\equiv x.$ Instead of 
a single universal premise $B_{\forall}$ 
we may have a finite conjunction of such 
premises. Instead of one space $P$ and one space $K$ we may have tuples of 
(potentially different) such spaces. 
\end{theorem} 
\noindent 
The main features of Theorem \ref{meta-corollary} are the following: 
\begin{itemize}
\item 
The extractability of a computable bound on both the premise as well as the 
conclusion (of course in practice these bounds will be different but by 
taking their maximum one always can obtain a common bound which makes things 
easier to state). In any concrete case, the 
bound extractable will not only be computable but 
of (usually low) 
subrecursive complexity depending on the principles used in the proof at 
hand. 
In our case we will obtain a rather simple bound in the end. 
\item 
The bound is highly uniform as it does not depend on $y\in K$ and on $x,T$ and 
$(X,d,W)$ altogether only via an upper bound $b$ on $d(z,T(z))$ and the 
distance of (approximate) fixed points of $T$ from $z$ (plus $\eta$ in the 
case of ${\cal A}^{\omega}[X,d,W,\eta]_{-b}$).
\item 
The assumption that $T$ has a fixed point is replaced by the existence of 
approximate fixed points (in some ball around $z$). The latter is usually more 
elementary to verify than the former and does not require the completeness of $X$ or closedness of $C$ (see Section \ref{main-results}).
\end{itemize}
\noindent 
The following theorem (essentially based on 
\cite{Rhoades(94),Schu(91),Schu(91B)}) 
is 
proved in \cite[Corollary 8]{K+Lambov} for the case of 
uniformly convex Banach spaces but its proof can be generalized to 
uniformly convex hyperbolic spaces.

\begin{theorem} \label{original-theorem}
Let $(X,d,W)$ be a (nonempty) uniformly convex hyperbolic space having a 
monotone modulus of uniform convexity $\eta$, $C\subseteq X$ be a nonempty convex subset and $T:C\rightarrow C$  be asymptotically nonexpansive with sequence $(k_n)\in [0,\infty)^{\N}$ 
where $\sum\limits^{\infty}_{i=0} k_i<\infty$. 
 Let $(\lambda_n)$ be a sequence in $[a,b]$ for $0<a<b<1.$

If  $T$ has a fixed point, then $T$ is $\lambda_n$-asymptotically regular.
\end{theorem} 

Since any convex subset of a hyperbolic space again is a hyperbolic space, it suffices to consider only functions $T:X\to X$. Then Theorem \ref{original-theorem} can be formalized as follows 
\[ \ba{l} \forall K,L,k\in\N\,\forall (\lambda_n)\in [0,1]^{\N}\,\forall 
(k_n)\in [0,K]^{\N}\,\forall x\in X\,\forall T:X\to X \\ 
\left( \ba{l} Mon(\eta)\,\wedge \, \forall n\in\N\,
\forall y,z\in X\big(d_X(T^ny,T^nz)\le_{\R} (1+k_n)d_X(y,z)\big) \\
\si \,\forall n\in\N\ \left(\ds\sum^n_{i=0} k_i\le K\right)\, \si\, L\ge 2\,\wedge\, \forall n\in\N\ \left(\ds\frac{1}{L} \le_{\R} 
\lambda_n\le_{\R}
1-\frac{1}{L}\right)\\ 
\wedge\, Fix(T)\not=\emptyset \\ 
\rightarrow \exists n\in\N
\forall m\in\N\ \big( d_X(x_{n+m},Tx_{n+m}\big)\le_{\R} 2^{-k}
\big)\ea \right). \ea \]
Here, $Mon(\eta)$ is the $\forall$-formula from \cite{L-2006} expressing that $\eta$ is monotone in the first argument (viewed as a rational number).
Then
\[\ba{l}   Mon(\eta)\,\wedge \, \forall n\in\N\,
\forall y,z\in X\big(d_X(T^ny,T^nz)\le_{\R} (1+k_n)d_X(y,z)\big) \\
\si \,\forall n\in\N\ \left(\ds\sum^n_{i=0} k_i\le K\right)\, \si\, L\ge 2\,\wedge\, \forall n\in\N\ \left(\ds\frac{1}{L} \le_{\R} 
\lambda_n\le_{\R}
1-\frac{1}{L}\right)
\ea 
\] is a finite conjunction of $\forall$-formulas and $[0,1]^{\N},[0,K]^{\N}$ 
are compact metric (and hence Polish) spaces ($\N$ is also covered by 
quantification over $P$ as mentioned above).

\begin{remark} Strictly speaking, $[0,K]^{\N}$ is not a single compact 
metric space but a sequence of such spaces as $K$ varies over $\N.$ However, 
this simple extension is also covered by (the proof of) Theorem 
\ref{meta-corollary}. 
\end{remark} \noindent 
The asymptotic nonexpansivity of $T$  
\[ \forall n\in \N\,\forall y,z\in X\,\big(d_X(T^ny,T^nz)\le_{\R} 
(1+k_n)d_X(y,z)\big) \] 
implies that $T$ is $(1+k_1)$-Lipschitz continuous. Since $k_1\le K,$ in fact 
$T$ is $(1+K)$-Lipschitz. So we do not need to add a Lipschitz constant as an 
extra input in order to be able to apply the logical metatheorem. \\[1mm]
Unfortunately, the conclusion  
 \[ \exists n\in\N \forall m\in\N \big( d_X(x_{n+m},Tx_{n+m})
\big)\le_{\R} 2^{-k}
\big) \] 
is not an $\exists$-formula, but only its weakened form 
\[ (*) \ \exists n\in\N\ \big( d_X(x_n,Tx_n)<_{\R} 2^{-k}
\big) \] is one. \\[1mm] 
Suppose now that the proof of Theorem \ref{original-theorem} 
can be formalized in 
${\cal A}^{\omega}[X,d,W,\eta]_{-b}$ (as is the case). 
Then the logical metatheorem stated 
above guarantees the extractability of a computable bound
$\Phi(K,L,b,\eta,k)$ such that the following holds in all (nonempty) uniformly 
convex hyperbolic spaces $(X,d,W,\eta)$ with monotone modulus $\eta$:

for all $K,L,k,b\in\N,(\lambda_n)\in [0,1]^{\N},(k_n)\in [0,K]^{\N}, 
x\in X,T:X\to X$ if $T$ is asymptotically nonexpansive with 
sequence $(k_n)$, $\lambda_n\in \left[\frac1L,1-\frac1L\right]$ for all $n\in\N$,  $\ds \sum_{k=0}^\infty k_n \leq K$ and 
\[ \forall \varepsilon >0 \, 
\big( Fix_{\varepsilon}(T,x,b)\not=\emptyset\big) \ \wedge \ 
d(x,Tx)\le b\] then 
\[ \exists n\le \Phi(K,L,b,\eta,k) \ \big( d(x_n,Tx_n)<2^{-k}\big). \]
The original convergence statement 
\[ (1) \quad \forall k\in\N\,\exists n\in\N\,\forall m\in\N \,
\big( d(x_{n+m},Tx_{n+m})<2^{-k}\big) \] 
can be rewritten as 
\[ (2) \quad \forall k\in\N\,\exists n\in\N\, \forall m\in\N \,
\forall i\in [n,n+m] \,\big( 
d(x_i,Tx_i)<2^{-k}\big). \]
$(2)$ clearly implies the so-called Herbrand normal form $(2)^H$ 
of $(2)$ 
\[ (2)^H \quad \forall k\in\N\,\forall g:\N\to\N\, \exists n\in\N\, 
\forall i \in [n,n+g(n)] \,\big( 
d(x_i,Tx_i)<2^{-k}\big). \] 
Ineffectively, also the converse is true, i.e. $(2)^H$ implies $(2)$ 
(and so also $(1)$): assume that $(2^H)$ is true. If $(2)$ would be false, 
then for some $k\in\N$
\[ \forall n\in\N\,\exists m_n\in\N \,
\exists i \in [n,n+m_n]\, (d(x_i,Tx_i)\ge 2^{-k}). \] 
Define $g(n):=m_n.$ Then $(2^H)$ applied to $g$ leads to a contradiction. 
Due to the ineffectivity of this argument, a bound on `$\exists n\in\N$' in 
$(2^H)$ cannot be 
converted effectively into a bound on `$\exists n\in\N$' in $(2).$ 
\[ \forall i\in [n,n+g(n)]\ \big( 
d(x_i,Tx_i)<2^{-k}\big) \] 
is equivalent to an $\exists$-formula (using that $<$ between real numbers is 
an existential formula and the universal quantifier over $i$ is bounded). 
Moreover, quantification over $\N^{\N}$ is covered (as mentioned above) even 
without any extra representation of the Baire space $\N^{\N}$ as a Polish 
metric space. \\[1mm] 
Hence one can apply the logical metatheorem also to the conclusion $(2^H)$ 
rather than just the special case 
\[ \exists n\in\N \ \big( d(x_n,Tx_n)<2^{-k} \big) \] 
which corresponds to $g(n)\equiv 0.$ As a result we can extract a computable 
bound $\Phi$ on `$\exists n\in\N$' in $(2^H)$ 
which in addition to $K,L,b,\eta,k$ also depends on $g,$ i.e. 
\[ (3) \ \exists n\le \Phi(K,L,b,\eta,k,g) \ \forall i\in [n,n+g(n)]\ 
 \big( d(x_i,Tx_i)<2^{-k}\big) \]
for all $g:\N\to\N.$ \\[2mm]
The rest of this paper is concerned with the construction of such a 
bound $\Phi(K,L,b,\eta,k,g)$, that is with the proof of Theorem \ref{Herbrand-ass-ne}.
We will carry out this construction directly by generalizing the 
reasoning from \cite{K+Lambov} rather than first proving Theorem \ref{original-theorem} and then extracting the bound from the proof. Note, however, that 
\cite{K+Lambov} was developed using the extraction algorithm underlying 
the proof of (earlier versions of) 
Theorem \ref{meta-corollary} (in its version for 
uniformly convex normed spaces). As 
$(3)$ implies $(2^H)$ and so (ineffectively) $(1)$ we will obtain as a 
corollary Theorem \ref{original-theorem}.
\begin{remark} At the time the paper \cite{K+Lambov} was written, the only 
logical metatheorems available (\cite{Kohlenbach(metapaper)}) required the 
boundedness of the convex subset in question. Only in 
\cite{Gerhardy/Kohlenbach3} the fact that the results in \cite{K+Lambov} 
did not require any global boundedness assumption could be accounted for by 
general logical theorems. In \cite{L-2006} this treatment was adapted to 
uniformly convex hyperbolic spaces, i.e. the context of the present paper.
\end{remark}

\section{Some technical lemmas}

In the following, $(X,d,W)$ is a hyperbolic space,  $C\se X$ a nonempty convex subset of $X$, $T:C\to C$ an asymptotically  nonexpansive mapping with sequence $(k_n)$, $(\lambda_n)$ is a sequence in $[0,1]$ and $(x_n)$  is the Krasnoselski-Mann iteration starting with $x\in C$.

\blem\label{KM-lemma}
Let $n\in\N, p\in C, \alpha>0, \gamma\geq \max\{d(Tp,p),d(T^np,p)\}$ and $K\geq k_m$ for all $m\in\N$. Then
\be
\item $d(T^nx_{n+1},x_{n+1})\leq  (1+k_n)d(T^nx_n,x_n)$ and 
\bea
d(x_{n+1},Tx_{n+1})\leq d(T^{n+1}x_{n+1},x_{n+1})+(1+K)^2d(T^nx_n,x_n).\label{Schu-hyperbolic}
\eea
\item \label{lemma-n-n+1} Assume that  for both $i=n$ and $i=n+1$ we have that
\[d(x_i,p)<\alpha\,\, \text{or}\,\,d(T^ix_i,x_i)< \alpha.\]
Then
\[d(x_{n+1},Tx_{n+1})< (1+(1+K)^2(2+K))\alpha+(1+K^2)\gamma.\]
\ee
\elem
\begin{proof}
\be
\item
\bua
d(T^nx_{n+1},x_{n+1})&= & d(T^nx_{n+1},(1-\lambda_n)x_n \oplus\lambda_n T^nx_n)\\
&\leq & (1-\lambda_n)d(T^nx_{n+1},x_n)+\lambda_n d(T^nx_{n+1},T^nx_n), \quad \text{by (W1)}\\
&\leq & (1-\lambda_n)d(T^nx_{n+1},T^nx_n)+(1-\lambda_n)d(T^nx_n,x_n)+\\
&& +\lambda_n d(T^nx_{n+1},T^nx_n)\\
&=& d(T^nx_{n+1},T^nx_n) +(1-\lambda_n)d(T^nx_n,x_n)\\
&\leq & (1+k_n)d(x_{n+1},x_n)+(1-\lambda_n)d(T^nx_n,x_n)\\
&=& (1+k_n)\lambda_n d(T^nx_n,x_n)+(1-\lambda_n)d(T^nx_n,x_n), \quad \text{by (\ref{prop-xylambda})}\\
&\leq & (1+k_n)d(T^nx_n,x_n).\\
d(x_{n+1},Tx_{n+1})&\leq &d(x_{n+1},T^{n+1}x_{n+1})+d(T^{n+1}x_{n+1},Tx_{n+1})\\
&\leq & d(x_{n+1},T^{n+1}x_{n+1})+(1+k_1)d(T^nx_{n+1},x_{n+1})\\
&\leq & d(x_{n+1},T^{n+1}x_{n+1})+(1+k_1)(1+k_n)d(T^nx_n,x_n)\\
&\leq & d(x_{n+1},T^{n+1}x_{n+1})+(1+K)^2d(T^nx_n,x_n).
\eua
\item
We have the following cases:
\be
\item[1.] $d(x_{n+1},p)< \alpha$. Then
\bua
d(x_{n+1},Tx_{n+1})&\leq & d(x_{n+1},p)+d(p,Tp)+d(Tx_{n+1},Tp) \\
&\leq & d(x_{n+1},p)+(1+k_1)d(x_{n+1},p)+d(p,Tp)\\
&=& (2+k_1)d(x_{n+1},p)+d(Tp,p)< (2+K)\alpha+\gamma.
\eua
\item[2.] $d(x_{n+1},p)\geq\alpha$. Then we must have $d(T^{n+1}x_{n+1},x_{n+1})<\alpha$. We distinguish two situations:
\be
\item[(a)] $d(T^nx_n,x_n)< \alpha$. Then, by (\ref{Schu-hyperbolic})
\bua
d(x_{n+1},Tx_{n+1})&\leq & d(T^{n+1}x_{n+1},x_{n+1})+(1+K)^2d(T^nx_n,x_n)\\
&< &  (1+(1+K)^2)\alpha.
\eua
\item[(b)] $d(x_n,p)<\alpha$. Then
\bua
d(T^nx_n,x_n)&\leq & d(T^nx_n,T^np)+d(T^np,p)+d(x_n,p)\\
&\leq & (2+k_n)d(x_n,p)+d(T^np,p)< (2+K)\alpha+\gamma.
\eua
Hence, using again (\ref{Schu-hyperbolic}),
\bua
d(x_{n+1},Tx_{n+1})&\leq & d(T^{n+1}x_{n+1},x_{n+1})+(1+K)^2d(T^nx_n,x_n)\\
&< & \alpha+(1+K)^2((2+K)\alpha+\gamma)\\
&\leq & (1+(1+K)^2(2+K))\alpha+(1+K^2)\gamma.
\eua
\ee
\ee
\ee
\end{proof}

\blem\label{main-technical-lemma}
Let $(X,d,W)$ be a uniformly convex hyperbolic space with a monotone modulus of uniform convexity $\eta$. Let $x,p\in C$ and $K\geq k_m$ for all $m\in\N$. Assume that $n\in\N, \alpha,\beta,\beta^*,\tilde{\beta},\gamma,\nu>0$ are such that
\bua
d(T^np,p)<\nu \leq 1, & \alpha\leq d(x_n,p)\leq \beta, \tilde{\beta},\beta^* \text{~and~} \alpha\leq d(x_n,T^nx_n).
\eua
Then 
\beq
d(x_{n+1},p) < d(x_n,p)+k_n\beta^*+\nu -2\alpha\lambda_n(1-\lambda_n)\eta\left((1+K)\tilde{\beta}+1,\frac{\alpha}{(1+K)\beta+1}\right). \label{eta-general}
\eeq
If, moreover, $\eta$ can be written as $\eta(r,\varepsilon)=\varepsilon\cdot\tilde{\eta}(r,\varepsilon)$ such that $\tilde{\eta}$ increases with $\varepsilon$ (for a fixed $r$), then
\beq
d(x_{n+1},p) < d(x_n,p)+k_n\beta^*+\nu -2\alpha\lambda_n(1-\lambda_n)\tilde{\eta}\left((1+K)\tilde{\beta}+1,\frac{\alpha}{(1+K)\beta+1}\right). \label{eta-tilde} 
\eeq
\elem
\begin{proof} 
Let $\ds r:=(1+k_n)d(x_n,p)+d(T^np,p), \,\,\eps:=\frac{\alpha}{(1+K)\beta+1} $ and $\ds\psi:=\frac{\alpha}{r}$. By hypothesis, 
$r< (1+K)\beta+1$, hence $0<\eps<\psi \leq 1$.

We note that
\begin{eqnarray*}
d(T^nx_n,p)&\leq & d(T^nx_n,T^np)+d(T^np,p)\leq r,\\
d(x_n,p)&\leq & r\\
d(x_n,T^nx_n) &\geq & \alpha = r\psi\geq r\eps.
\end {eqnarray*}
We get that
\bua
\rho(x_{n+1},p)&= & \rho((1-\lambda_n)x_n \oplus\lambda_n T^nx_n,p) \\
&\leq & \big(1-2\lambda_n(1-\lambda_n)\eta(r,\eps)\big)\cdot r, \quad\text{by Lemma \ref{eta-prop-1}.\ref{uc-ineq-Groetsch}}\\
&\leq & \big(1-2\lambda_n(1-\lambda_n)\eta((1+K)\tilde{\beta}+1,\eps)\big)\cdot r, \\
&&\quad\text{since~~}r< (1+K)\tilde{\beta}+1\text{~~and~~}\eta\text{~~is monotone}\\
&= & r - 2r\lambda_n(1-\lambda_n)\eta((1+K)\tilde{\beta}+1,\eps)\\
&\leq & r - 2\alpha\lambda_n(1-\lambda_n)\eta((1+K)\tilde{\beta}+1,\eps), \quad \text{since~} r\geq\alpha\\
&<& d(x_n,p)+k_n\beta^*+\nu -2\alpha\lambda_n(1-\lambda_n)\eta((1+K)\tilde{\beta}+1,\eps).
\eua

Assume now that $\eta(r,\varepsilon)=\varepsilon\cdot\tilde{\eta}(r,\varepsilon)$ and $\tilde{\eta}$ increases with $\varepsilon$. 
Applying again Lemma \ref{eta-prop-1}.\ref{uc-ineq-Groetsch} and the monotonicity of $\eta$, but with $\psi$ instead of $\eps$, we obtain
\bua
d(x_{n+1},p)&\leq & \big(1-2\lambda_n(1-\lambda_n)\eta((1+K)\tilde{\beta}+1,\psi)\big)\cdot r
\\
&=& \big(1-2\lambda_n(1-\lambda_n)\psi\tilde{\eta}((1+K)\tilde{\beta}+1,\psi)\big)\cdot r\\
&= & r - 2\alpha\lambda_n(1-\lambda_n)\tilde{\eta}((1+K)\tilde{\beta}+1,\psi)\\
&\leq & r -2\alpha\lambda_n(1-\lambda_n)\tilde{\eta}((1+K)\tilde{\beta}+1,\eps),\quad \text{since~} \eps<\psi\\
&< & d(x_n,p)+k_n\beta^*+\nu -2\alpha\lambda_n(1-\lambda_n)\tilde{\eta}((1+K)\tilde{\beta}+1,\eps).
\eua
\end{proof}

We shall use also the following quantitative lemmas on sequences of real numbers.

\blem\label{basic-lemma}
Let  $(a_n)_{n\geq 0}$ be a real sequence. Then
\beq
\forall\eps>0 \forall g:\N\to\N\left(a_{g^M(0)}\geq 0\ra \exists i<M\big(a_{g^i(0)}-a_{g^{i+1}(0)}\leq \eps)\right),
\eeq
where $ \ds M :=\left\lceil \frac{a_0}{\varepsilon} \right\rceil. $
As a consequence, 
\beq
\forall\eps>0 \forall g:\N\to\N\left(\forall n\leq \Theta(a_0,\eps,g)(a_n\geq 0)\ra \exists N\leq \Theta\big(a_N-a_{g(N)}\leq \eps)\right),
\eeq
where $\Theta(a_0,\eps,g):= \ds\max\{g^i(0): i\leq M\}$. Moreover,  $N=g^i(0)$ for some $i<M$.
\elem
\begin{proof}
Let $\eps >0, g:\N\to\N$ be such that $a_{g^M(0)}\geq 0$.  Assume by contradiction that  $a_{g^i(0)}-a_{g^{i+1}(0)}>\eps$ for all $i\in \ol{0, M-1}$.  By adding these inequalities, we get that $\ds a_0-a_{g^M(0)}>M\eps=\left\lceil \frac{a_0}{\eps}\right\rceil\cdot\eps\geq  \frac{a_0}{\eps}\cdot\eps=a_0$,
hence $a_0-a_{g^M(0)}>a_0$, which is a contradiction, since $a_{g^M(0)}\geq 0$. 
\end{proof}

The following lemma is a special case of 
\cite[Lemma 17]{K+Lambov}.\footnote{Corrections to \cite{K+Lambov}: 1) 
In Lemma 15 and below $\lfloor\cdot\rfloor$ should be 
$\lceil\cdot\rceil.$ P.164, line 4: `$j\le m+g(m)$' should be `$j\le m+g(m)+1$' 
and, consequently, in Theorem 22 and Corollary 28 `$h=\lambda n.(g(n+1)+n+1$' 
should be `$h=\lambda n.(g(n+1)+n+2$' and in Corollary 25 `$n\le\Phi_1$' must 
be replaced by `$n\le 2\Phi_1$'.} 

\bprop\label{prop-quant-qihou-2-seq-v2}
Let $A_1,A_2\ge 1,B_1,B_2, C_1,C_2\geq 0$ 
and define for any $\theta>0$ and for any $g:\N\to \N$
\beq
\Psi(A_1,A_2, B_1,B_2, C_1,C_2, g, \theta) := h^M(0) , \label{quant-qihou-Phi-2-seq-v2}
\eeq
where 
\[\ba{l}
\ds h(n) :=  g(n)+n, \quad D_i :=  (A_i+C_i)\exp(B_i), \\[0.1cm]
\ds M :=  \ds\left\lceil\frac{3(4B_1D_1 + 4C_1 + D_1+4B_2D_2 + 4C_2 + D_2)}{\theta}\right\rceil.
\ea\]
Let $(a_n), (b_n), (c_n),(\alpha_n),(\beta_n),(\gamma_n)$ be real sequences  such that for all $n\leq \Psi$,
\[ a_n, b_n, c_n,\alpha_n,\beta_n,\gamma_n\geq 0,\quad 
a_{n+1} \le (1+b_{n})a_{n}+c_{n}\quad \text{and}\quad \alpha_{n+1} \le (1+\beta_{n})\alpha_{n}+\gamma_{n}\]
and, moreover, 
\[a_{0} \le A_1 \,\,, \,\, \alpha_0\leq A_2,\,\,\,\,  \ds\sum_{n=0}^\Psi b_{n} \le B_1\,\,,\,\,\sum_{n=0}^\Psi \beta_{n} \le B_2,\,\,\,\, \ds\sum_{n=0}^\Psi c_{n} \le C_1\,\,,\,\,  \sum_{n=0}^\Psi \gamma_{n} \le C_2.\]
Then the following holds:
\be
\item $a_n\leq D_1, \alpha_n\leq D_2$ for all $n\leq \Psi+1$;
\item for all $\theta\in(0,1]$ and all $g:\N\to\N$,
\bua
 \exists N\le \Psi \forall i,j\in[N,N+g(N)]\left( |a_j - a_i| \le \theta \,\,\si \,\, |\alpha_j - \alpha_i| \le \theta\right).\label{quant-qihou-2-v2}
\eua
Moreover, $N=h^i(0)$ for some $i<M$.
\ee
\eprop

\section{Proof of Theorem \ref{Herbrand-ass-ne}}

Let $\eps\in(0,1],g:\N\to\N$ be arbitrary and $K,L,x\in C, b, h:\N\to\N, M, D, \theta, f(K),\Phi$ as given in the hypotheses of Theorem \ref{Herbrand-ass-ne}. 
Let us remark that $\ds\frac{\eps}{f(K)((1+K)D+1)}<\frac{1}{2}<1$ and, moreover, $\ds\theta\leq\frac{\eps}{L^2f(K)}<1$.

Since  $x\in C$ and $b>0$ satisfy (\ref{hyp-app-b-x}), there exists $p\in C$ such that
\beq
d(x,p)\leq b \quad \si \quad d(p,Tp)\leq \frac{1}{2^{\Phi }(\Phi +K)}.\label{bound-Tp-p}
\eeq
Since 
\bua
d(T^np,p) &\leq & d(T^np,T^{n-1}p)+ d(T^{n-1}p,p)\leq (1+k_{n-1})d(p,Tp)+d(T^{n-1}p,p),
\eua
it follows that for all $1\leq n \leq\Phi $, 
\bua
d(T^np,p) &\leq & \sum_{i=0}^{n-1}(1+k_i)d(p,Tp)=d(p,Tp)\left(n+\sum_{i=0}^{n-1}k_i\right)\\
&\leq & (n+K)\cdot\frac{1}{2^{\Phi }(\Phi +K)}\leq \frac{1}{2^{\Phi }} \leq\frac1{2^n}.
\eua
Let us consider the sequences:
\bua
a_n:=d(x_n,p), \quad \alpha_0:=KD+2,\, \alpha_n:=KD+2-\sum_{i=0}^{n-1}\left(k_i D+\frac1{2^i}\right) \text{~~for~~} n\geq 1. 
\eua
Then for all $n\leq \Phi$,  we have that $0\leq \alpha_{n+1}\leq \alpha_n$ and
\bua
0\leq a_{n+1}&=& d((1-\lambda_n)x_n \oplus\lambda_n T^nx_n,p)\leq (1-\lambda_n)d(x_n,p)+\lambda_n d(T^nx_n,p)\\
&\leq & (1-\lambda_n)d(x_n,p)+\lambda_n d(T^nx_n,T^np)+\lambda_n d(T^np,p)\\
&\leq &(1-\lambda_n)d(x_n,p)+\lambda_n(1+k_n)d(x_n,p)+\lambda_n d(T^np,p)\\
&\leq&(1+k_n)d(x_n,p)+d(T^np,p)\leq (1+k_n)a_n+\frac1{2^n}.
\eua
It is easy to verify that we can apply Proposition \ref{prop-quant-qihou-2-seq-v2} with $a_n,\alpha_n$ given as above, $\ds b_n:=k_n, c_n:=\frac1{2^n}, \beta_n:=\gamma_n:=0$, $A_1:=b,B_1:=K, C_1:=2, A_2:=KD+2,   B_2:=C_2:=0$,  $\tilde{g}(n):=g(n+1)+2$ and $\theta,\Phi$ as above. 

It follows by Proposition \ref{prop-quant-qihou-2-seq-v2} that 
\beq
a_n\leq (A_1+C_1)\exp(B_1)=D \text{~~for all~~} n\leq\Phi \label{an-leq-D}
\eeq
and that there exists $N_0\leq \Phi, N_0=h^s(0)$ for some $s<M$ such that
\beq
\forall i, j\in[N_0,N_0+g(N_0+1)+2]\left(|a_j - a_i|\leq \theta \,\si\, |\alpha_j-\alpha_i|\leq \theta\right).\label{basis-Psi}
\eeq
In fact, since the sequence $(h^n(0))$ is strictly increasing, we have that $N_0=h^s(0)<h^M(0)=\Phi$, so $N_0+1\leq \Phi$. 

Let $N:=N_0+1$. In the following, we shall prove that $N$ satisfies (\ref{Herbrand-ass-ne-conclusion}), that is
\bua
\forall m\in [N, N+g(N)] \left (d(x_m,Tx_m) <\eps \right).
\eua
Let $m\in [N,N+g(N)]$. Then $m-1,m,m+1\in [N_0,N_0+g(N_0+1)+2]$, so
we can apply (\ref{basis-Psi}) with $i\in\{m-1,m\}$ and $j=i+1$  to get that
\beq
|d(x_{i+1},p)-d(x_{i},p)|=|a_{i+1}-a_{i}|\leq\theta \,\,\text{~and~} \,\,  k_{i}D+\frac1{2^i}=|\alpha_{i+1}-\alpha_{i}|\leq \theta
\label{quant-m-m-1}
\eeq
Moreover, 
\bua
m-1< m\leq N_0+1+g(N_0+1)< h(N_0)=h^{s+1}(0)\leq h^M(0)\leq\Phi ,\label{m-leq-Phi}
\eua
Let $i\in\{m-1,m\}$ and assume that $\ds d(x_i,p)\geq \frac{\eps}{f(K)}$ and $\ds d(T^ix_i,x_i)\geq \frac{\eps}{f(K)}$. Then  
\bua
d(T^ip,p) \leq\frac1{2^\Phi}< \frac1{2^i} \leq 1,   \quad \frac{\eps}{f(K)}\leq d(T^ix_i,x_i), \quad
 \frac{\eps}{f(K)}\leq d(x_i,p)\leq D \text{~~(by (\ref{an-leq-D}))},
\eua
so we can apply Lemma \ref{main-technical-lemma}, (\ref{eta-general}) with $\ds\alpha:=\frac{\eps}{f(K)},\, \nu:=\frac1{2^i},\, \beta:=\beta^*:=\tilde{\beta}:=D$, the definition of $\theta$ and the fact that  $\ds\lambda_i(1-\lambda_i)\geq \frac1{L^2}$ to get that
\bua
d(x_{i+1},p)& < &   d(x_i,p)+k_iD+\frac1{2^i} -2\theta.
\eua
It follows that
\bua
2\theta <  d(x_i,p)-d(x_{i+1},p)+k_iD+\frac1{2^i}=a_i-a_{i+1}+k_iD+\frac1{2^i}\leq 2\theta, \text{~ by (\ref{quant-m-m-1})},
\eua
that is a contradiction.

Hence, for both $i=m$ and $i=m-1$.
\bua
d(x_i,p)< \frac{\eps}{f(K)}\,\, \,  \text{~or~} \, \,\,  d(T^ix_i,x_i) < \frac{\eps}{f(K)}.\label{both-i-i+1}
\eua 

Finally, applying Lemma \ref{KM-lemma}.\ref{lemma-n-n+1} with $\ds n:=m-1, \alpha:=\frac{\eps}{f(K)},\gamma:=\frac1{2^\Phi } $ it follows that
\bua
d(x_{m},Tx_{m})&< &  (1+(1+K)^2(2+K))\frac{\eps}{f(K)}+(1+K^2)\frac1{2^\Phi }\\
&=& \frac{\eps}{2}+(1+K)^2\frac1{2^\Phi }, \text{~~by the definition of~} f(K)\\
&<& \frac{\eps}{2}+(1+K^2)\frac1{2^m} \text{~~since~~} m<\Phi \\
&\leq & \frac{\eps}{2}+(1+K)^2\theta, \text{~~ since~} \frac1{2^m}\leq\theta, \text{~by (\ref{quant-m-m-1})}\\
&< & \eps.
\eua
since $\ds (1+K)^2\theta\leq \frac{(1+K)^2\eps}{L^2f(K)}< \frac{\eps}{2}$.
\qed


\begin{thebibliography}{}

\bibitem{Borwein+Reich+Shafrir-92}
J. Borwein, S. Reich, I. Shafrir, Krasnoselski-Mann iterations in normed 
spaces. Canad. Math. Bull 35 (1992),  21-28.

\bibitem{Bridson+Haefliger-book}
M. Bridson, A. Haefliger, Metric spaces of non-positive curvature.
Grundlehren der Mathematischen Wissenschaften, 319. Springer-Verlag, Berlin, 1999. xxii+643 pp. 

\bibitem{Clarkson-1936}
J.A. Clarkson, Uniformly convex spaces.
 Trans. Amer. Math. Soc.  40  (1936),  no. 3, 396--414. 


\bibitem{Gerhardy/Kohlenbach3} 
P. Gerhardy, U. Kohlenbach,  General logical metatheorems for functional 
analysis. Trans. Amer. Math. Soc. 360 (2008), no. 5,  2615-2660. 

\bibitem{Goebel(02)} K. Goebel, Concise course on fixed point theory. Yokohama Publishers, Yokohama, 2002. iv+182 pp. 


\bibitem{Goebel+Kirk-72} 
K. Goebel, W.A. Kirk, A fixed point theorem for asymptotically nonexpansive 
mappings. Proc. Amer. Math. Soc.  35  (1972), 171--174. 


\bibitem{Goebel+Kirk-1983}
K. Goebel, W.A. Kirk, Iteration processes for nonexpansive mappings. In:  
S. P. Singh, S. Thomeier, B. Watson (eds.), Topological methods in nonlinear functional analysis (Toronto, 1982),  115--123, Contemp. Math., 21, Amer. Math. Soc., Providence, RI, 1983. 

\bibitem{Goebel+Reich-book}
K. Goebel, S. Reich, Uniform convexity, hyperbolic geometry, and 
nonexpansive mappings. Monographs and Textbooks in Pure and Applied Mathematics, 83. Marcel Dekker, Inc., New York, 1984. ix+170 pp.

\bibitem{Khamsi-2003}
M.A. Khamsi, On asymptotically nonexpansive mappings in hyperconvex metric 
spaces.   Proc. Amer. Math. Soc.  132  (2004),  no. 2, 365--373. 

\bibitem{Kirk-1982}
W.A. Kirk, Krasnosel'skii iteration process in hyperbolic spaces,
Numer. Funct. Anal. and Optimiz. 4 (1982), 371-381.

\bibitem{Kirk-2004}
W.A. Kirk, Geodesic geometry and fixed point theory II. In:   J. Garcia Falset, 
E. Llorens Fuster, B. Sims (eds.), International Conference on Fixed Point Theory and Applications (Valencia, 2003), 113--142, Yokohama Publ., Yokohama, 2004.

\bibitem{Kohlenbach(A)} U. Kohlenbach, Analysing proofs in
analysis. In: W. Hodges, M. Hyland, C. Steinhorn, J. Truss,
(eds.),  Logic: from foundations to applications (Staffordshire, 1993),  225--260, Oxford Sci. Publ., Oxford Univ. Press, New York, 1996.


\bibitem{Kohlenbach(metapaper)} U. Kohlenbach, Some logical metatheorems 
with applications in functional analysis. Trans. Amer. Math. Soc. 357 
(2005),  no. 1, 89-128.

\bibitem{Kohlenbach(06a)} U. Kohlenbach, Effective uniform bounds from proofs 
in abstract functional analysis.  In: B. Cooper, B. Loewe, A. Sorbi, (eds.), New Computational Paradigms: Changing Conceptions of What is Computable, 223-258, Springer-Verlag, Berlin, 2008.

\bibitem{Kohlenbach(book)} U. Kohlenbach,  Applied Proof Theory: Proof 
Interpretations and their Use in Mathematics. 
Springer Monographs in Mathematics. Springer-Verlag, Berlin, 2008. xix+532pp.

\bibitem{K+Lambov}
U. Kohlenbach, B. Lambov, Bounds on iterations of asymptotically 
quasi-nonexpansive mappings. In: J. Garcia Falset, E. Llorens Fuster, B. Sims 
(eds.),  International Conference on Fixed Point Theory and Applications (Valencia, 2003),  143--172, Yokohama Publ., Yokohama, 2004.

\bibitem{L-07-JMAA}
L. Leu\c stean, A quadratic rate of asymptotic regularity for CAT(0)-spaces.
J. Math. Anal. Appl.  325  (2007),  no. 1, 386--399. 


\bibitem{L-2006} L. Leu\c{s}tean, Proof mining in $\R$-trees and hyperbolic 
spaces.  In: G. Mints and R. de Queiroz (eds.), Proceedings of the 13th Workshop on Logic, Language, Information and Computation (WoLLIC 2006) (Stanford, 2006), 95--106, Electron. Notes Theor. Comput. Sci., 165, Elsevier, Amsterdam, 2006.

\bibitem{Qihou(01A)} L. Qihou, Iteration sequences for asymptotically 
quasi-nonexpansive mappings with error member. J. Math. Anal. Appl. 
259 (2001), 18-24. 

\bibitem{Qihou(02)} L. Qihou, Iteration sequences for asymptotically 
quasi-nonexpansive mapping with an error member of uniform convex
Banach space. J. Math. Anal. Appl. 266 (2002), 468-471. 

\bibitem{Reich+Shafrir-1990}
S. Reich, I. Shafrir, Nonexpansive iterations in hyperbolic spaces. 
Nonlinear Analysis 15 (1990), 537-558.

\bibitem{Rhoades(94)} B.E. Rhoades, Fixed point iterations for 
certain nonlinear mappings. J. Math. Anal. Appl. 183 (1994), 118-120. 

\bibitem{Schu(91)} J. Schu, Iterative construction of fixed points of 
asymptotically nonexpansive mappings. J. Math. Anal. Appl. 158 (1991), 407-413.  

\bibitem{Schu(91B)} J. Schu, Weak and strong convergence to fixed 
points of asymptotically nonexpansive mappings. Bull. Austral. Math. 
Soc. 43 (1991), 153-159. 

\bibitem{Spector(62)} C. Spector, Provably recursive functionals of 
analysis: a consistency proof of analysis by an extension of principles 
formulated in current intuitionistic mathematics. In:  J.C.E. Dekker (ed.), 1962 Proc. Sympos. Pure Math., Vol. V, 1--27, Amer. Math. Soc., Providence, R.I., 1962.

\bibitem{Takahashi-70}
W. Takahashi, A convexity in metric space and nonexpansive mappings, I. 
Kodai Math. Sem. Rep. 22 (1970), 142--149.

\bibitem{Tao(07)} T. Tao, Soft analysis, hard analysis, and the finite 
convergence principle. Essay posted May 23, 2007. Available at: 
http://terrytao.wordpress.com/2007/05/23/soft-analysis-hard-analysis-and-the-finite-convergence-principle/.

\bibitem{Tao(07a)} T. Tao, Norm convergence of multiple ergodic averages 
for commuting transformations. arXiv:0707.1117v1 [math.DS] (2007). To appear in Ergodic Theory and Dynamical Systems. 

\end{thebibliography}
\end{document}